\documentclass[copyright]{eptcs}

\usepackage{amssymb,amsmath,amsthm}
\usepackage{tikz}
\usetikzlibrary{cd}

\usepackage{graphicx}
\usepackage{cleveref}
\usepackage{listings}
\usepackage{enumitem}

\usepackage{iftex}

\ifpdf
  \usepackage{underscore}         
  \usepackage[T1]{fontenc}        
\else
  \usepackage{breakurl}           
\fi

\usepackage{framed}

\newcounter{instruction}
\newenvironment{instruction}[1][]
{
    \refstepcounter{instruction}
    \begin{framed}
    \noindent\textbf{Instruction \theinstruction\ifx&#1&\else: #1\fi}
    \vspace{0.5em}
    \newline
}
{
    \end{framed}
}

\title{The Educational Proof Assistant Waterproof in an Introductory Proof Course: Proof Construction and Learning Processes}
\author{Pim Otte
\institute{Utrecht University\\ Utrecht, The Netherlands}
\institute{Eindhoven University of Technology\\ Eindhoven, The Netherlands}
\email{p.j.otte@uu.nl}
\and
Rogier Bos
\institute{Utrecht University\\ Utrecht, The Netherlands}
\email{r.d.bos@uu.nl}
\and
Johan Commelin
\institute{Utrecht University\\ Utrecht, The Netherlands}
\email{j.m.commelin@uu.nl}
\and
Jim Portegies
\institute{Eindhoven University of Technology\\ Eindhoven, The Netherlands}
\email{j.w.portegies@tue.nl}
}

\providecommand{\publicationstatus}{}

\def\titlerunning{The Educational Proof Assistant Waterproof in an Introductory Proof Course}
\def\authorrunning{P. Otte, R. Bos, J. Commelin and J. Portegies}
\begin{document}
\maketitle

\begin{abstract}
  We study the use of an educational proof assistant in an introductory proof course through a quasi-experiment in a varied setting:
  multiple teachers, students with different study programs, and a mixed Dutch-English language environment.
  First-year university students are known to struggle with writing proofs.
  Waterproof is a proof assistant that is designed to support the transfer of skills to paper proofs by
  working with controlled natural language.
  We focus on the students' ability to construct valid mathematical proofs, and on their learning process.
  We study this through in-class observation, surveys, and analysis of student performance and proof structure.
  We present evidence that effects of using an educational proof assistant carry over to the pen-and-paper context,
  even when the assistant is English and the proof is given in Dutch.
  We also present evidence that suggests students in the Mathematics-Computer Science program achieve higher grades
  when using Waterproof.
  Our most important conclusion is that an educational proof assistant can help students be more explicit in their proofs.
  As students self-selected into using Waterproof rather than being randomly assigned, these results are suggestive rather than causal.
\end{abstract}

\section{Introduction}\label{sec:introduction}

In this work, we study the use of an educational proof assistant in education
in a first-year undergraduate introductory proofs course.
For the purposes of this work, an educational proof assistant is a proof assistant adapted 
for education. We chose to evaluate \emph{Waterproof}\footnote{\url{https://impermeable.github.io}} in an introductory proof course at
Utrecht University.
\emph{Waterproof} was conceived and is still being developed at Eindhoven University of Technology~\cite{Wemmenhove2024waterproof}.

\subsection{Motivation}

Students have many difficulties with learning to produce (formal) mathematical proofs, as documented by e.g. Moore~\cite{moore1994making}.
He finds the major sources to be concept understanding, mathematical language and notation, and getting started on a proof.
Proof assistants provide the user with feedback on their proofs, but skills obtained using a proof assistant do not
necessarily transfer to writing proofs on paper (see Knobelsdorf~et~al.~\cite{knobelsdorf2017theorem}).
This inspires the adaptation of proof assistants to be more suitable for teaching mathematical proof writing.

\subsection{Relation to existing research}

This section situates this work in existing research along various lines. 
We describe the relations with existing tools and research.

The use of proof assistants in teaching goes back to the development and use of \emph{Mizar-PC} for teaching logic by Matuszewski, 
Rudnicki and Trybulec in 1975~\cite{matuszewski2005mizar}. Following this, many more tools have been developed and 
evaluated. For an extensive overview, we refer to Tran Minh~et~al.~\cite{tranminh2025proof}.
We discuss selected examples in relation to the current work.

Beyond the above application of Mizar in teaching logic, there are many more domain-specific applications
of educational use of proof assistants. 
For logic, this includes 
the use of \emph{Isabelle} by Villadsen and Jacobsen~\cite{villadsen2021using}, 
the \emph{ProofWeb} system by Kaliszyk~et~al.~\cite{wiedijk2007teaching}, and
the \emph{TryLogic} tutorial by Terrematte and Marcos~\cite{terrematte2015trylogictutorialapproachlearning}.
In the case of geometry, we highlight work on automated theorem proving
using \emph{Geogebra} (Botana~et~al.~\cite{botana2015automated}, Quaresma and Santos~\cite{quaresma2019computer}), 
\emph{GeoProof} by Narboux~\cite{narboux2007graphical}.
All of these differ from Waterproof mainly in targeting a specific area: Waterproof concerns mathematical
proof in the general setting.

There are also various existing approaches with respect to user interfaces and input methods for
educational proof assistants.
Tools like \emph{D$\exists\forall$DUCTION} by Le Roux~(see \cite{kerjean2022utilisation}) and
\emph{Actema} by Donato~et~al.~\cite{donato2023integrating} function using \emph{proof by pointing},
inviting mouse-based interactions.
\emph{Paperproof} by Karunus~and~Kovsharov~\cite{paperproof} visualizes Lean proofs along with changes in the proof state.
\emph{Diproche} by Carl~et~al.~\cite{carl2022natural}, and \emph{Verbose Lean} by Massot~\cite{massot2024teaching} are most
like Waterproof in that all of these get input through controlled natural language.
In these works, the tools specifically use German and French to accommodate the students,
whereas in this study we apply an English tool in a mixed language setting: materials are in English and lectures (largely) in Dutch.

With respect to prior evaluations of educational proof assistants in the classroom, we build on
two studies that both analyze the proof structure of students' proofs.
The impact of Waterproof on students' proof structure was evaluated by Wemmenhove~et~al.~\cite{wemmenhove2026comparative}.
Thoma~and~Iannone~\cite{thoma2022learning} explored the characteristics of proof writing by first-year undergraduate students who participated in Lean workshops.
We build on these works by analyzing a proof that was given by students after the period in which Waterproof was used with most intensity.
The main difference between these works and the current work is that we have a control group.

\subsection{Research Questions}

In this work we consider the following two research questions:

\begin{enumerate}
    \item[RQ1] How does the use of an educational proof assistant in an introductory proof course affect students' ability to construct high-quality mathematical proofs?
    \item[RQ2] How does the use of an educational proof assistant in an introductory proof course affect the students'  learning process?
\end{enumerate}

The first research question concerns the issue whether the tool actually helps. 
As a proxy measure for the quality of proof, we use students' performance in the course 
and the structure of their proofs.

The second question captures the idea of knock-on effects.  
This includes various factors that could explain a positive or negative effect with respect to their proof writing ability.
We choose to zoom in on three aspects of the students' learning process: the perception of feedback they receive, their time-on-task,
and their help-seeking behavior. Since the software provides feedback, we hypothesize students using it
will perceive receiving more and more frequent feedback, than the control group that receives only feedback from the teacher and teaching assistants.
Anecdotal classroom observations suggest that proof assistants can be very engaging,
which supports a hypothesis that students will spend more time on task when using an educational proof assistant.
With respect to help-seeking behavior, there are reasons to believe both more and fewer questions will arise; 
the proof assistant might answer some questions, but its use might also raise questions.

\subsection{Overview of results}

We find that students who actively use Waterproof tend to use more explicit proof steps. 
This is supported by the combination of grades and analyses of the quiz.
We do not reach a conclusion with respect to changes to the students' learning process.
There are minor indications students do not experience errors from Waterproof as
feedback and that students using Waterproof ask slightly more questions than students who do not. 
We also find evidence that students in the Mathematics/Computer Science program obtain higher grades
after practicing with Waterproof.

\section{Methods}

\subsection{Context}

The context for this research is the first-year course ``Bewijzen in de Wiskunde'', taught at Utrecht University.
This is an introductory proof course which teaches students the basics of mathematical proof, along
with basic set theory, functions, relations, modular arithmetic and limits.
This course is taught in seven parallel groups, by different teachers. Two of the authors, Commelin and Otte, were among the teachers.
The teachers collaborated on course materials and assessment, but taught separate mixed lecture-tutorial classes.
These classes took place twice a week, with a duration of 3 hours and 45 minutes.
The students taking the course are first-year undergraduate students, with a fair variety among their study programs. 
The majority of students study mathematics, and a large fraction of them intend to double major in physics, computer science or economics,
and some of them study mathematics and applications specifically.

\subsection{Intervention}\label{subsec:intervention}

Commelin gave an introduction to Waterproof in the first lesson, showing the functionality.
Otte did the same, and additionally presented selected proofs live in Waterproof during the first weeks of the course. 
Commelin and Otte provided students with both examples and exercises in
Waterproof.\footnote{Available here: \url{https://github.com/impermeable/introduction-to-proof-sheets/tree/9ab77fa6a01411ab2e9f79d1958f56683521bcfc}}
Use of Waterproof by students was voluntary.

Waterproof is an educational proof assistant based on the Rocq Prover. Students complete exercises through 
the use of controlled natural language. Waterproof allows asking for help and will give feedback to
enforce a particular style of proof. We provide a visual impression in \Cref{app:screenshot}.

In all groups, we asked students to participate in the research. We asked them to fill in two surveys, and 
 for consent for being observed, as well as the use of their proofs and grades for analysis.
The protocol was approved by the Science-Geosciences Ethics Review Board of Utrecht University.

In summary, this setup provides a
quasi-experiment\footnote{A quasi-experiment is an experiment, except participants are not randomly assigned to the treatment or control group.} in a realistic environment.
This setup splits the students into three groups: a full control group from the 5 groups without Waterproof,
a relative control group of students who were in Commelin's and Otte's groups, but did not choose to use it, and
an intervention group of students who actively used Waterproof. This allows a view on both the effect of
the intervention among different teachers, as well as insight on the students who chose or chose not to use Waterproof.

\subsection{Quiz}\label{subsec:quiz}

The course was assessed through two in-class quizzes, two homework exercises and a final written exam.
In this research, we focus on the performance of students in the first quiz.
This quiz was given in the middle of the second week of teaching, and aligns best
with active use of Waterproof by students, as described in \Cref{subsec:engagement}.

We analyzed the grades from this quiz with a statistical analysis using a linear mixed model.
All grades reported in this study are on the Dutch scale from 1 to 10, where a grade of 5.5 or higher is a pass.
We control for prior knowledge and skills in two ways: We use a pretest to control for prior knowledge
of the learning goals tested in the quiz, and the secondary school exam grade\footnote{For students from the Dutch secondary system this is the grade for the course ``Mathematics~B''; students from other systems reported an equivalent grade on the 1--10 scale.} to control for general mathematical ability.
The pretest was an exercise similar in form and content to the first quiz, administered at the start of the course before the
relevant material had been taught, in order to establish students' (lack of) prior knowledge of the tested subject matter.
These are both included as effects in the linear mixed model. The teacher and study program of students are included as
random effects. The main independent variable is the number of exercises completed with Waterproof
and the dependent variable is the first quiz grade. All numerical variables are standardized.

We analyzed the proof structure in a similar manner as used by Wemmenhove~et~al.~\cite{wemmenhove2026comparative}.
We transcribed the proofs and labeled each step with codes corresponding
to logical rules or proof steps. We coded proofs line by line, and labeled
a line as `continued' if it continued a rule or step from the previous line.
Incorrect or unclear steps were not labeled at all and are reported as `no code'.
See \Cref{app:quiz-coding} for a fully coded reference solution, and the full coding scheme.

\subsection{Surveys}

Participants filled in a survey at the start and end of the course.
The survey at the start asked about their mathematics performance in secondary education,
their study program, gender and age. The full entry survey is attached as \Cref{app:entry-survey}.
The final survey focused on various elements. We surveyed the time spent on the course, 
perception of feedback given, reasons for continuing or stopping to use Waterproof. 
The full final survey is attached as \Cref{app:final-survey}.

We aggregated the survey responses and used them in fitting the Linear Mixed Model described in \Cref{subsec:quiz} 
in our answer to RQ1.

We use the survey in various ways to support our answer to RQ2.
The final survey contained multiple questions about the quantity of feedback (Q 4aci), frequency (Q 4dfh) and whether students
got feedback when needed (Q 4beg). The responses were collected on a five-point Likert scale. We aggregate these responses per
category. We aggregate the hours spent on the course from Q6 across active Waterproof users, inactive Waterproof users
and the control group. Due to participants' differing interpretation of whether the question
was asking time spent per week or in total, we multiplied all answers below 20 hours by 8.

\subsection{Observation}

In-class observations were performed by two observers. They observed
the class taught by Otte, and one other where Waterproof was not used.
These classes were selected because the respective teachers have roughly the same amount of experience teaching the course,
as well as teaching in general.
The final 2 hours and 45 minutes of each class were observed, since the first hour was
generally lecture, often with a break at the end and consistent between the two groups.
The observers each selected 4 students among participants of the research at the beginning of
the observation period, aiming to balance the amount of times each student was observed. 
They observed these students at least every 5 minutes and annotated
their observations with codes. The observation plan and codes are available in \Cref{app:observationplan}.
They also recorded all questions asked by the students chosen for observation.

We analyze the observations to support our answer to RQ2.
We group the observations to lead to time-series data for an observed session per student.
We categorized the observations
by students spending time on the course, mathematics unrelated to the course, unrelated tasks, or ambiguous.
The results are averaged per 2 hour 45 minutes observation period, and grouped
by whether Waterproof was used in this period. This means observed sessions from the same student on different
days might be in different categories, depending on whether they used Waterproof in this time.

\section{Results}\label{sec:results}

\subsection{Participant population}

We describe the participant population.
Of the 199 students enrolled in the course, 80 consented to participate in the research and completed the entry survey.
Of them, 52 completed the exit survey.
We use the population of students who completed the first survey where this is sufficient,
and the population of students who completed the exit survey where we need the full data.
Among the participant population, 46 identified as male,
32 identified as female, 1 as neither and 1 chose not to disclose. The average age was 18 years and 10 months.
This is in line with expectations for the general population of first-year mathematics students at Utrecht University.
Because the groups are not randomized, there is a biased distribution of students among the groups with respect to their study program.

\subsubsection{Engagement with Waterproof}\label{subsec:engagement}

Of the 34 students in the Waterproof groups, 21 completed one or more exercises in Waterproof,
and 16 completed at least five. We consider the latter group ``active Waterproof users'',
and the 18 other students as ``inactive Waterproof users''. The latter group were in classrooms
where others were using Waterproof, but (barely) engaged with it themselves.
Engagement was best until the third week of the course. In the fourth week, no Waterproof exercises
were available.
Very few students continued using Waterproof after this, even when there were exercises available.
\Cref{tab:participants-by-programme} provides an overview of the number of students in the different subgroups.

\begin{table}[htbp]
\centering
\caption{Participants by study programme and Waterproof engagement.}
\label{tab:participants-by-programme}
\begin{tabular}{lrrr}
\hline
Study Programme & Active Waterproof & Inactive Waterproof & Control \\
\hline
Mathematics (no double bachelor) &  4 &  3 & 14 \\
Mathematics and applications &  0 &  0 &  9 \\
Mathematics/Computer Science &  4 &  9 &  0 \\
Mathematics/Economics &  1 &  0 &  5 \\
Mathematics/Physics &  5 &  4 & 14 \\
Other: &  2 &  2 &  4 \\
Total & 16 & 18 & 46 \\
\hline
\end{tabular}
\end{table}

\subsection{RQ1: Ability of students to construct proofs}

\subsubsection{Grades of students in the first quiz}

\Cref{tab:quiz1-by-programme} presents the students' average grades of the first quiz, split out by 
study programme and Waterproof engagement.

\begin{table}[htbp]
\centering
\caption{Mean Quiz 1 grades by study programme and Waterproof engagement}
\label{tab:quiz1-by-programme}
\begin{tabular}{lrrr}
\hline
Study Programme & Active Waterproof & Inactive Waterproof & Control \\
\hline
Mathematics (no double bachelor) & 6.12 & 5.75 & 5.96 \\
Mathematics and applications &  &  & 4.61 \\
Mathematics/Computer Science & 8.00 & 6.78 &  \\
Mathematics/Economics & 5.50 &  & 6.65 \\
Mathematics/Physics & 7.20 & 4.12 & 6.75 \\
Other: & 5.75 & 4.50 & 5.62 \\
\hline
\end{tabular}
\end{table}
\Cref{tab:fixed-effects-lmm} presents the estimated effect, standard error, the degrees of freedom, the $t$-statistic, and the $p$-value
of the variables in the linear mixed model.
We observe that the secondary school exam grade has a statistically significant ($p \approx .001$) large (Cohen's d $\approx 0.844$) positive effect on the first quiz grade.
The number of exercises done in Waterproof has a positive but not statistically significant effect ($p \approx 0.1$).

\begin{table}[htbp]
\centering
\caption{Fixed effects of the Linear Mixed Model (all variables normalized)}
\label{tab:fixed-effects-lmm}
\begin{tabular}{lrrrrr}
\hline
Term & Estimate & Std. Error & df & $t$ & $p$ \\
\hline
(Intercept) & 0.001 & 0.188 &  5.3 & 0.008 & 0.994 \\
Number of Waterproof exercises & 0.189 & 0.113 & 57.1 & 1.669 & 0.101 \\
Secondary school exam grade & 0.373 & 0.103 & 73.5 & 3.617 & 0.001 \\
Pretest grade & 0.034 & 0.108 & 73.9 & 0.319 & 0.751 \\
\hline
\end{tabular}
\end{table}

\subsubsection{Proof structure in the first quiz}

\Cref{tab:coding-full-by-group} shows the most interesting results of the proof analysis.
Each code corresponds to a logical rule or proof step. The codes relevant here are
$\exists$-intro (proving an existential by providing a witness, e.g.\ ``Choose $y := 14$''),
$\vee$-elim (splitting a disjunction into cases),
signpost-case (announcing or opening a case in a case split, e.g.\ ``Case~1: $\neg p(x)$''), and
signpost-concl (stating an intermediate or final conclusion, e.g.\ ``In both cases \ldots'').
We highlight the higher prevalence of $\exists$-intro, $\vee$-elim and signpost-case among active Waterproof users,
as well as fewer signpost-concl steps.
These codes are relevant to the experiment because they capture how explicitly students structure their proofs:
signposting cases and introducing concrete witnesses are exactly the steps that Waterproof's controlled natural language requires,
whereas a signpost-concl step can leave the underlying reasoning implicit.
For examples of these steps we refer to \Cref{app:quiz-coding}, and the full coding scheme is available in \Cref{app:complete-coding}.

\begin{table}[htbp]
\centering
\caption{Average occurrences of the most interesting proof codes per participant, by group}
\label{tab:coding-full-by-group}
\begin{tabular}{lrrr}
\hline
Code & Active & Inactive & Control \\
\hline
$\exists$-intro & 2.31 & 1.41 & 1.48 \\
$\vee$-elim & 1.44 & 0.65 & 0.70 \\
signpost-case & 4.31 & 3.12 & 2.78 \\
signpost-concl & 1.62 & 2.06 & 2.20 \\
\hline
\end{tabular}
\end{table}

\subsection{Proof structure across language}

During the grading of the written final exam, we observed that some Waterproof users still adhered to language and structure that
was clearly Waterproof-inspired, though in Dutch. 
We present one of the clearest cases of this in \Cref{tab:dutchwp}. A translation is provided in \Cref{app:translation}. 
In particular we note the strict domain check in the existential introduction.

\begin{table}[h]
\centering
\begin{tabular}{|l|}
\hline
  Te bewijzen is dat $\forall y \in f(A), \exists x \in A, y = g(x)$. \\
  Neem $y \in f(A)$. \\
  $f(A) = \{f(x)\mid x \in A\}$, dus $\exists x' \in A, y = f(x)$. \\
  Verkrijg deze $x' \in A$. \\
  Er geldt dus $y \in f(x')$. \\
  We moeten bewijzen dat $\exists x \in A, y = g(x)$. \\
  Kies $x := x'$. \\
  (Inderdaad, $x \in A$.) \\
  We moeten bewijzen dat $y = g(x')$. \\
  Per definitie van $g$ geldt dat $g(x')=f(x')$, dus het is waar dat $y = g(x')$. \\
  We concluderen dat $g$ surjectief is. \\
\hline

\end{tabular}
\caption{Dutch proof of surjectivity of a function given on the exam by an active Waterproof user, showing a Waterproof-inspired proof}
\label{tab:dutchwp}
\end{table}

\subsection{RQ2: Impact on learning process}

\subsubsection{Perceptions about feedback}

\begin{figure}[htbp]
  \centering
  \includegraphics[width=0.8\linewidth]{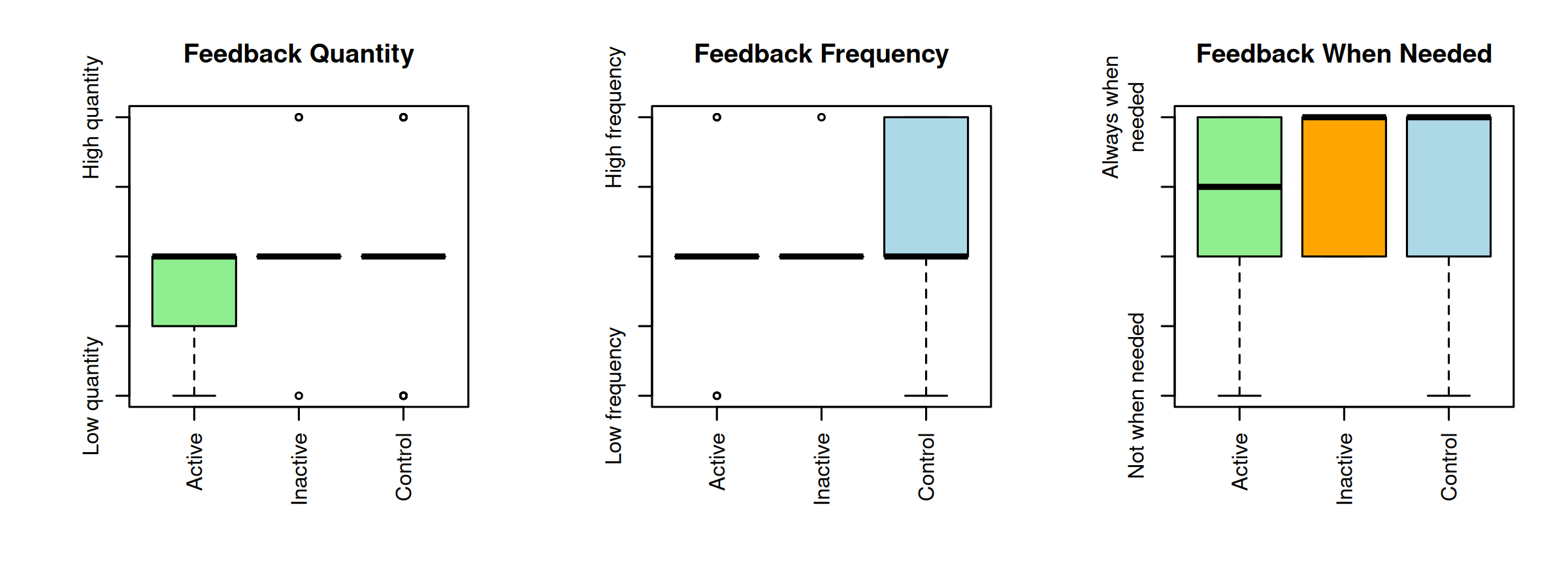}
  \caption{Survey results on perceptions about feedback}
  \label{fig:feedback}
\end{figure}

\Cref{fig:feedback} presents the answers to the questions about feedback, aggregated and normalized.
We note that active Waterproof users perceive a worse situation with respect to feedback;
we analyze this seemingly counter-intuitive observation further in \Cref{subsec:answer-rq2}.

\subsubsection{Time on task}

\begin{figure}[htbp]
  \centering
  \includegraphics[width=0.6\linewidth]{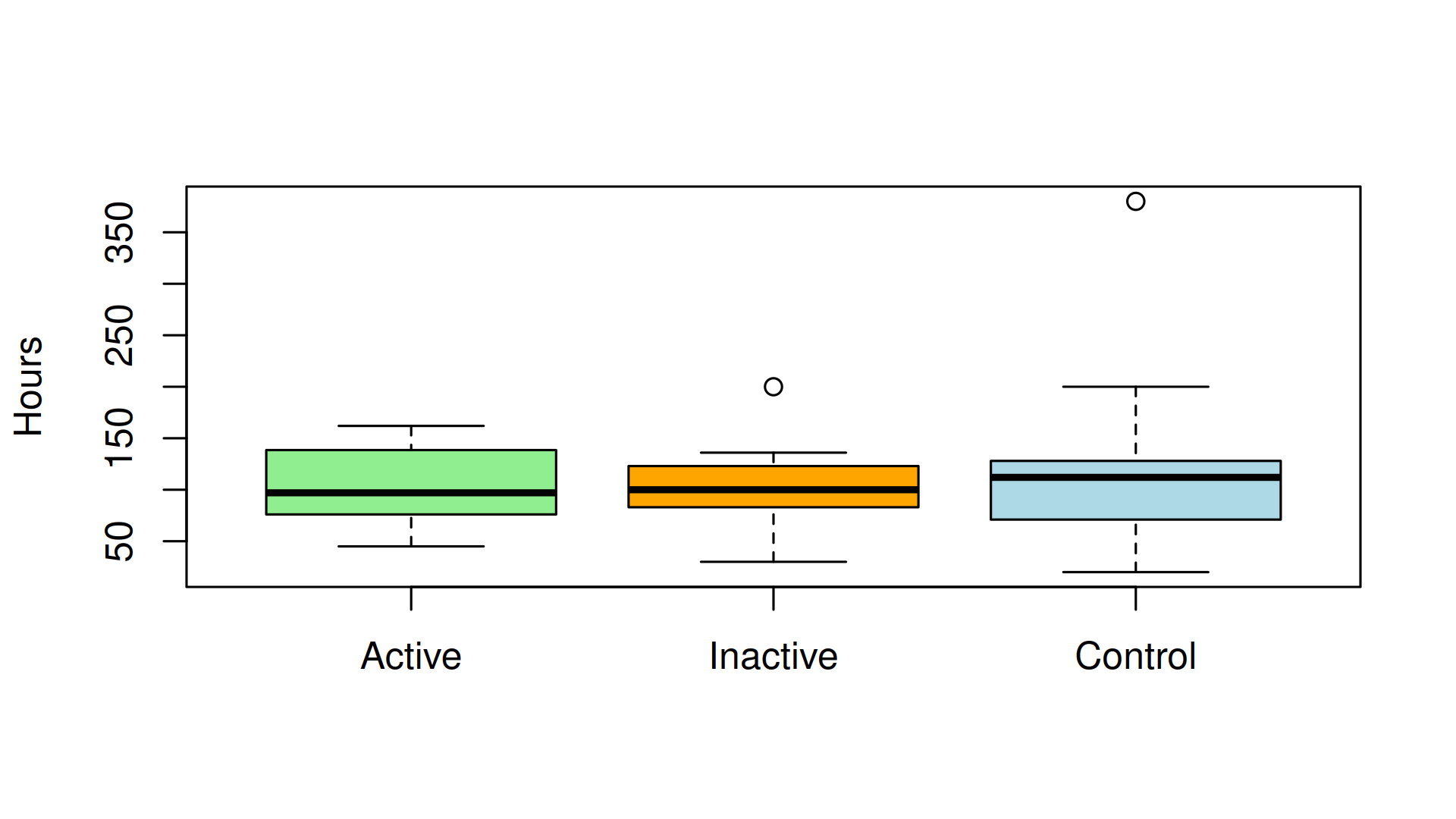}
  \caption{Survey results on time spent on the course}
  \label{fig:time}
\end{figure}

\Cref{fig:time} presents the self-reported time spent on the course, aggregated by Wateproof use.

\begin{figure}[htbp]
  \centering
  \includegraphics[width=0.6\linewidth]{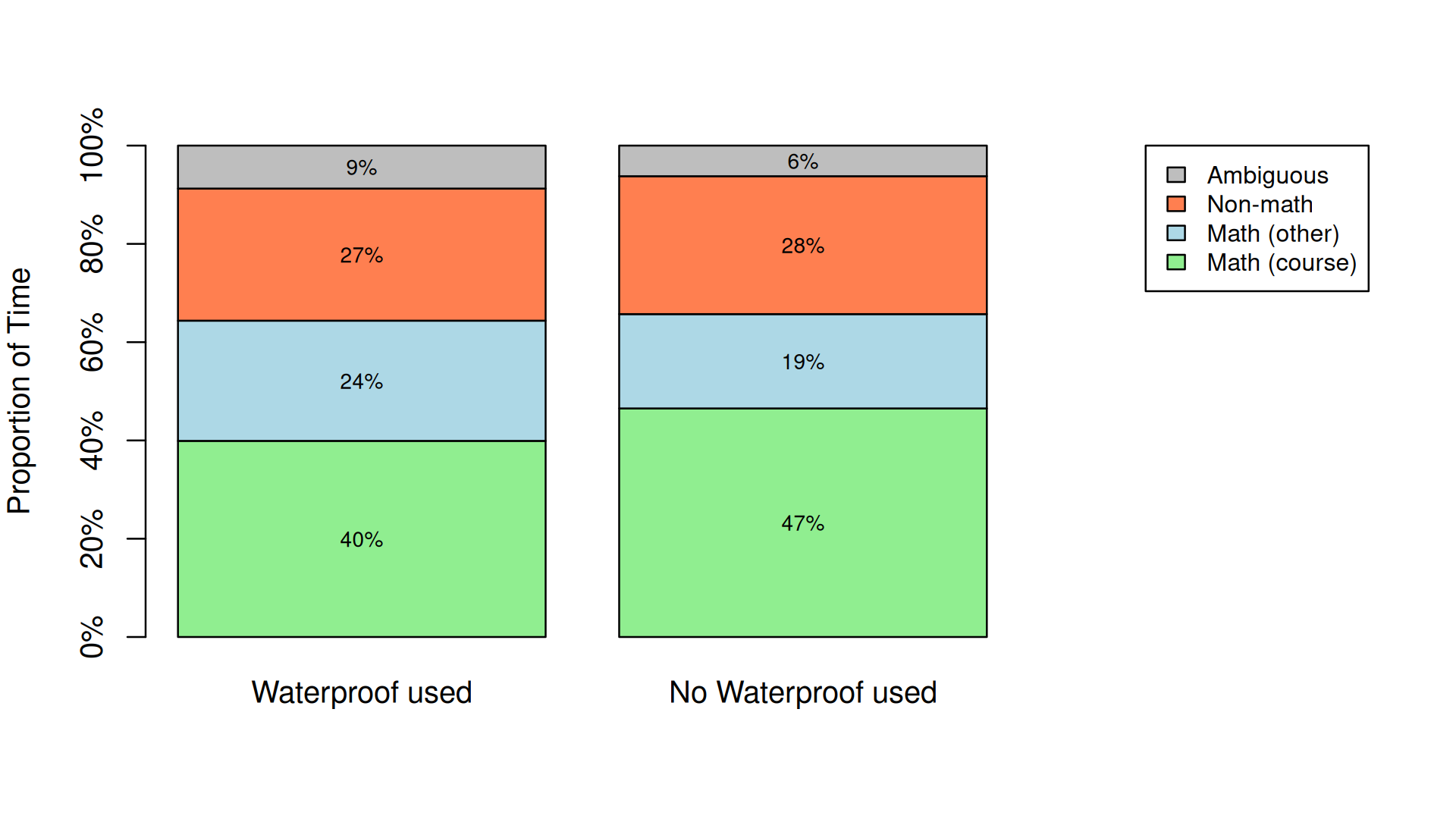}
  \caption{Observations on time spent in class. The vertical axis shows the proportion of observed in-class time spent
  in each activity category, averaged over the observation sessions and grouped by whether Waterproof was used during the session.}
  \label{fig:timeontask}
\end{figure}

For students' time spent in the classroom, we relied on observation.
The total number of sessions observed was 136, of which 22 included some amount of Waterproof use.
\Cref{fig:timeontask} shows, for each group, the proportion of observed in-class time that students spent on
the course mathematics, on other mathematics, on non-mathematical activities, or on activities we could not classify (ambiguous).

\subsubsection{Help-seeking behavior}

\begin{figure}[htbp]
  \centering
  \includegraphics[width=0.6\linewidth]{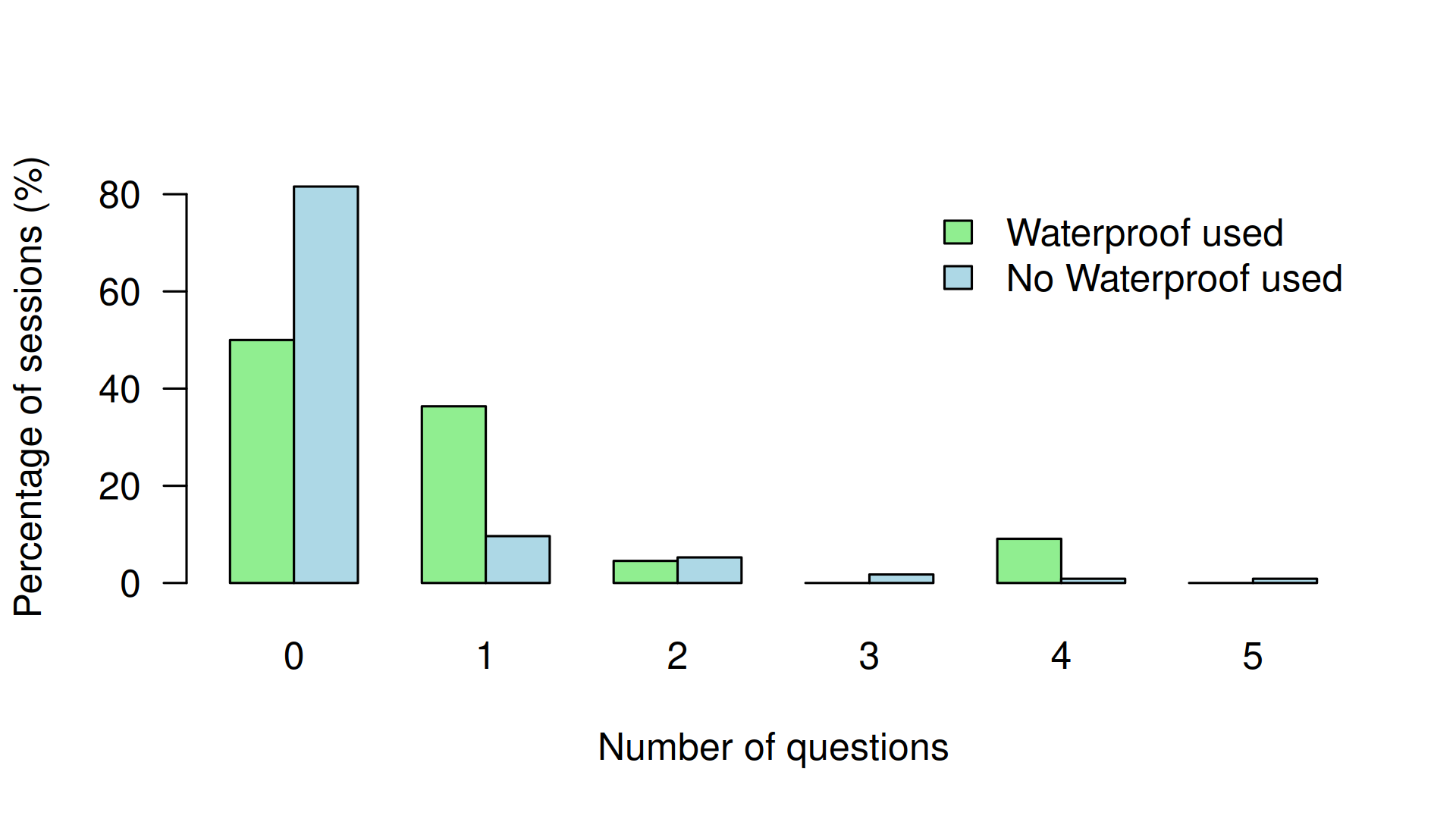}
  \caption{Observations on number of questions asked by students}
  \label{fig:questions}
\end{figure}

In \Cref{fig:questions}, we report the average number of questions asked by students per 2 hours 45 minutes observation period, grouped by whether
they used Waterproof in this period.
We note that Waterproof students asked more questions.

\section{Discussion}

We conclude by answering the research questions and placing our research in a broader context.

\subsection{Answer to RQ1: Students' proofs become more explicit}

We conclude that students who actively use Waterproof produce more explicit proofs than students who do not.
The fact that signpost-case, $\vee$-elim and $\exists$-intro steps occur more often (see \Cref{tab:coding-full-by-group})
 is a good sign: 
These steps are the core steps in the reference solution (see \Cref{app:quiz-coding}). If we combine this with the grades 
not being significantly higher, this points to the difference being something that is not in the rubric.
We also note that active Waterproof users make fewer (intermediate) conclusions
through signpost-concl steps (\Cref{tab:coding-full-by-group}).
We conclude that active Waterproof users more often explicitly stated their steps
rather than stating an intermediate conclusion and leaving the exact steps implicit.
We propose this is a consequence of Waterproof reinforcing patterns in making these steps. 
For example, signposting cases is obligatory in Waterproof.

We claim that finding some students transfer Waterproof syntax to Dutch shows clear evidence of 
transfer of skills obtained using Waterproof to the context of paper proofs. 
The proof in \Cref{tab:dutchwp} shows that the skills obtained in the Waterproof context have transcended
the language barrier in addition to transferring to paper. This suggests this is not just
imitation or an automatism, but really affected the way this student writes down proofs.

We cannot draw a definitive conclusion with respect to the performance of Mathematics-Computer Science
students, due to their absence in the control group. However, \Cref{tab:quiz1-by-programme} suggests that
they might benefit so much from Waterproof that their grades are positively affected.

Our findings in RQ1 are in line with earlier findings. Wemmenhove~et~al.~\cite{wemmenhove2026comparative}
find increased incidence of signposting. In their study, it mainly manifests as signposting
of the goal, rather than of cases. We propose this is due to the nature of the exercises studied in
both works. Thoma~and~Iannone~\cite{thoma2022learning} find
``the (often overt) breakdown of the proof goal'' as a characteristic of proofs by Lean users, which
matches the notion of making explicit steps.

\subsection{Answer to RQ2: Inconclusive impact on the learning process}\label{subsec:answer-rq2}

\Cref{fig:feedback} suggests that Waterproof users might not perceive 
Waterproof hints or error messages as feedback.
Generally, Waterproof produces a lot of errors, and if participants perceived these as feedback, 
they would have the perception of receiving more feedback, not less.
\Cref{fig:time} and \Cref{fig:timeontask} show no clear trends either way with respect to how students spend their time. 
\Cref{fig:questions} shows that students ask slightly more questions in classes where they
used Waterproof, but we think the more relevant aspect here is the low absolute amount of
questions from either group.

\subsection{Limitations}

A major limitation is the effect of different teachers. Otte and Commelin have both been strongly influenced
in their teaching by working with proof assistants, and this is hard to untangle from the effect of Waterproof.

While a quasi-experimental setup is an improvement over existing work, our implementation has a limitation to generalizability. 
Students are not randomly assigned to a group, but they are grouped partially by study program. Anecdotally, some students in a double degree program perform better than
students studying only mathematics.
We mitigate this in the statistical analysis by controlling with the pretest and secondary school exam grade.
We also remark that there is a dual role of teacher and researcher, which risks affecting 
both teacher and student behavior.

\section{Conclusion}

Our main finding is the increased explicitness in proofs by Waterproof users.
We stress that this and all our findings stem from a quasi-experiment in which students self-selected into using Waterproof
rather than being randomly assigned. The active, inactive and control groups may therefore differ systematically,
so our results should be read as evidence pointing towards an effect rather than as proof of a causal one.
The evidence provided with respect to improved performance of students of the mathematics/computer science program,
and increased number of questions in sessions with Waterproof use provides subjects for further research.
Our finding with respect to low absolute amount of questions asked by students provides reason to advise
the teachers to think about how to encourage students to ask more questions.

\section{Acknowledgements}

We thank Xander Baaijens and Michiel Kaptijn for performing the in-class observations and helping with the transcriptions.
We thank the teachers of ``Bewijzen in de Wiskunde'' for their cooperation.
We thank the many student developers who worked on Waterproof.
We thank Paige Randall North for her feedback on the writing.

\bibliographystyle{eptcs}
\bibliography{thedu26}

\pagebreak
\appendix

\section{Waterproof screenshot}\label{app:screenshot}

\begin{figure}[htbp]
  \centering
  \includegraphics[width=\linewidth]{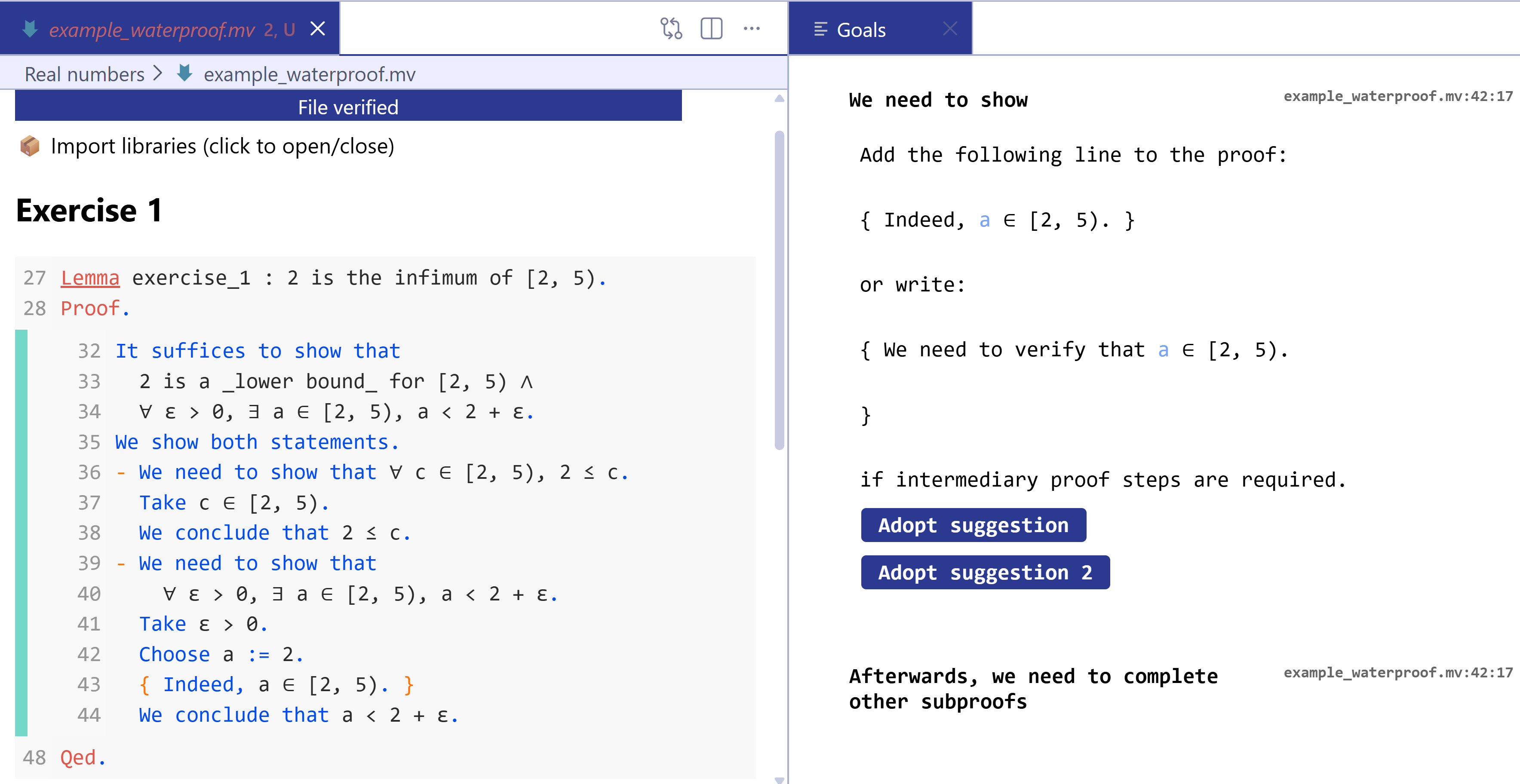}
  \caption{Screenshot showing a Waterproof proof and example of feedback.}
  \label{fig:screenshot}
\end{figure}

\section{Quiz and coding}\label{app:quiz-coding}

We provide the exercise contained in the studied quiz, along with the coded reference solution, in Waterproof syntax.

        \textbf{Exercise:} \\
       Let $p(x)$ and $q(x)$ be predicates with a free variable $x \in \mathbb{R}$. \\
        Prove that for all $x \in \mathbb{R}$, if
        $$p(x) \Rightarrow \bigl((\neg q(x) \Rightarrow q(3)) \vee q(14)\bigr)$$
        then
        $$\bigl(\neg p(x)\bigr) \vee \bigl(\exists y \in \mathbb{R}, q(y)\bigr)$$

        \textbf{Solution:} \\

        \begin{tabular}{|p{0.7\textwidth}|p{0.25\textwidth}|}
        \hline
        \textbf{Step} & \textbf{Code} \\
        \hline
        Take $x \in \mathbb{R}$. & $\forall$-intro \\
        \hline
        Assume that $p(x) \Rightarrow \bigl((\neg q(x) \Rightarrow q(3)) \vee q(14)\bigr)$. & $\Rightarrow$-intro \\
        \hline
        We need to show $\neg p(x) \vee \exists y \in \mathbb{R}, q(y)$. & signpost-goal \\
        \hline
        Either $\neg p(x)$ or $p(x)$. & use-thm, $\vee$-elim \\
        \hline
        \textbf{- Case} $\neg p(x)$. & signpost-case \\
        \hline
        \quad We conclude that $\neg p(x) \vee \exists y \in \mathbb{R}, q(y)$. & signpost-concl \\
        \hline
        \textbf{- Case} $p(x)$. & signpost-case \\
        \hline
        \quad It holds that $(\neg q(x) \Rightarrow q(3)) \vee q(14)$. & $\Rightarrow$-elim \\
        \hline
        \quad Either $q(14)$ or $(\neg q(x) \Rightarrow q(3))$. & $\vee$-elim \\
        \hline
        \quad \textbf{+ Case} $q(14)$. & signpost-case \\
        \hline
        \qquad It suffices to show $\exists y \in \mathbb{R}, q(y)$. & or-intro \\
        \hline
        \qquad Choose y := 14. & $\exists$-intro \\
        \hline
        \qquad \{ Indeed, $14 \in \mathbb{R}$\} & domain-check \\
        \hline
        \qquad We conclude that q(14). & signpost-concl \\
        \hline
        \quad \textbf{+ Case} $(\neg q(x) \Rightarrow q(3))$. & signpost-case \\
        \hline
        \qquad Either $q(x)$ or $\neg q(x)$ & use-thm, $\vee$-elim \\
        \hline
        \qquad\textbf{* Case} $q(x)$. & signpost-case \\
        \hline
        \quad\qquad It suffices to show $\exists y \in \mathbb{R}, q(y)$. & or-intro \\
        \hline
        \quad\qquad Choose y := x. & $\exists$-intro \\
        \hline
        \quad\qquad \{ Indeed, $x \in \mathbb{R}$\} & domain-check \\
        \hline
        \quad\qquad We conclude that q(x). & signpost-concl \\
        \hline
        \qquad\textbf{* Case} $\neg q(x)$. & signpost-case \\
        \hline
        \quad\qquad It holds that $(\neg q(x) \Rightarrow q(3))$. & $\Rightarrow$-elim \\
        \hline
        \quad\qquad It suffices to show $\exists y \in \mathbb{R}, q(y)$. & or-intro \\
        \hline
        \quad\qquad Choose y := 3. & $\exists$-intro \\
        \hline
        \quad\qquad \{ Indeed, $3 \in \mathbb{R}$\} & domain-check \\
        \hline
        \quad\qquad We conclude that q(3). & signpost-concl \\
        \hline
        \end{tabular}

The full coding scheme is presented in \Cref{tab:coding-examples}.

\begin{table}[htbp]
\centering
\caption{Proof coding taxonomy with example phrases}
\label{tab:coding-examples}
\small
\begin{tabular}{lp{5cm}p{6.5cm}}
\hline
\textbf{Code} & \textbf{Description} & \textbf{Example} \\
\hline
$\exists$-elim      & Eliminates an existential; names a witness                          & ``Let $a$ be such that $q(a)$'' \\
$\exists$-intro     & Proves an existential by providing a witness                        & ``Choose $y = 14$'' \\
$\forall$-elim      & Instantiates a universally quantified statement at a specific value & Applying $\forall x,\,q(x)$ at $x=3$ to get $q(3)$ \\
$\forall$-intro     & Introduces an arbitrary element to prove a universal                & ``Let $x\in\mathbb{R}$ be arbitrary'' \\
$\neg$-intro        & Derives a negation by assuming the positive and reaching a contradiction & ``Assume $q(x)$; \ldots contradiction; hence $\neg q(x)$'' \\
$\Rightarrow$-elim  & Applies modus ponens: from $A$ and $A\Rightarrow B$ derives $B$     & ``From $p(x)$ and $p(x)\Rightarrow(\ldots)$ it follows that $\ldots$'' \\
$\Rightarrow$-intro & Assumes the antecedent to prove an implication                      & ``Assume $p(x)\Rightarrow((\neg q(x)\Rightarrow q(3))\vee q(14))$'' \\
$\vee$-elim         & Splits on a disjunction into cases                                  & ``We split into two cases'' \\
$\vee$-intro        & Introduces a disjunction by establishing one disjunct               & ``Thus $\neg p(x)\vee\exists y\in\mathbb{R},\,q(y)$'' \\
domain-check        & Verifies that a value lies in the required domain                   & ``$14\in\mathbb{R}$'' \\
signpost-case       & Announces or opens a case in a case split                           & ``Case 1: $\neg p(x)$'' \\
signpost-concl      & Concludes, summarises, or announces end of proof                    & ``In both cases $\neg p(x)\vee\exists y\in\mathbb{R},\,q(y)$'' \\
signpost-goal       & States what remains to be proved                                    & ``It suffices to show $\exists y\in\mathbb{R},\,q(y)$'' \\
unfold              & Expands or rewrites a definition                                    & ``$A\Rightarrow B\equiv\neg A\vee B$'' \\
use-thm             & Applies a named theorem or rule without proof                       & ``$p(x)\vee\neg p(x)$ (law of excluded middle)'' or just ``We know that $p(x)\vee\neg p(x)$'' \\
continued           & Continuation of a proof step across multiple lines                  & A formula on its own line following the line coded with the rule \\
scrap               & Draft or scratch work                                               & Any line under a ``KLAD:'' heading \\
no code             & Organizational text, invalid steps                          & ``Bewijs:'', ``Z.O.Z.'', page-turn markers \\
\hline
\end{tabular}
\end{table}

\section{Complete coding results}\label{app:complete-coding}

\Cref{tab:coding-full-by-group-full} presents the full results of the coding of the proof from the first quiz.

\begin{table}[htbp]
\centering
\caption{Average occurrences of all proof codes per participant, by group}
\label{tab:coding-full-by-group-full}
\begin{tabular}{lrrr}
\hline
Code & Active & Inactive & Control \\
\hline
$\exists$-elim & 0.00 & 0.06 & 0.00 \\
$\exists$-intro & 2.31 & 1.41 & 1.48 \\
$\forall$-elim & 0.00 & 0.06 & 0.00 \\
$\forall$-intro & 0.19 & 0.12 & 0.22 \\
$\neg$-intro & 0.00 & 0.00 & 0.02 \\
$\Rightarrow$-elim & 1.19 & 0.88 & 0.74 \\
$\Rightarrow$-intro & 0.69 & 0.41 & 0.46 \\
$\vee$-elim & 1.44 & 0.65 & 0.70 \\
$\vee$-intro & 1.00 & 1.12 & 1.20 \\
domain-check & 0.94 & 0.76 & 0.57 \\
signpost-case & 4.31 & 3.12 & 2.78 \\
signpost-concl & 1.62 & 2.06 & 2.20 \\
signpost-goal & 1.38 & 1.29 & 1.13 \\
unfold & 0.00 & 0.18 & 0.17 \\
use-thm & 1.25 & 0.82 & 0.70 \\
continued & 2.94 & 2.71 & 4.48 \\
scrap & 0.50 & 1.06 & 1.59 \\
no code& 3.62 & 3.53 & 5.22 \\
\hline
\end{tabular}
\end{table}

\section{Observation plan}\label{app:observationplan}
These are the instructions for TAs observing students during "Bewijzen in de Wiskunde".
The observations focus on two main aspects: Time on task and help-seeking behavior.

The mixed lecture-tutorial format runs for 4 hours each session. The first hour is generally lecture+break, and is
unobserved. The last 3 hours are mostly tutorial mixed with some lecture, and are observed.
The first session is not observed as a whole, as it is needed to set up the research.

The TAs tasked with observation have observation as the primary task: They are allowed to answer questions,
and help students, but their main task is to observe and record data.

\begin{instruction}[Attendance and Setup]
    At the start of the observation period, one hour after the start of the lecture, 
    the TA records attendance of students for purposes of course administration.
    This time is also used to put the names of students that participate in the study 
    in a useful order for observations during the session.
    Select four students to observe during the lecture. Feel free to ask students
    to swap seats if it helps to be able to observe four of them in one go.
    Try to vary the students being observed over different classes.
\end{instruction}

\subsection{Time on Task}

\begin{instruction}[Recording Time on Task]
Repeatedly observe each student participating in the research and record whether they are actively engaged with course material.
Write down your observation in the provided Excel sheet. If possible, also add the codes immediately, but if the observation is clearly
written down, you can add the codes later.

Use the following codes:
\begin{itemize}
    \item \textbf{Ex}: The student is working on an exercise.
    \item \textbf{WP}: Student is working on Waterproof exercise.
    \item \textbf{Work}: The student working on the course in another way.
    \item \textbf{OC}: Student is working on some other course (specify which course in comment).
    \item \textbf{OT}: Student is doing something unrelated to any course (e.g. smartphone use, discussing weekend plans, sleeping).
    \item \textbf{A}: The student is absent from the classroom.
    \item \textbf{W}: The student is waiting on something.
    \item \textbf{O}: Other, explain in comment.
\end{itemize}

Record the observations in an Excel sheet in OneDrive with the name observations\_groupnumber\_YYYYMMDD.xlsx.
The columns should be: Student name, Time (rounded to nearest minute), Observation, Code, Comment.

Aim to continuously observe each of the four students, recording times to the nearest minute.
Short behavior does not need to be recorded, e.g. if a student checks their phone for a few seconds.

We will process the data assuming that codes hold until the next observation for that student.

\end{instruction}

\subsection{Help-Seeking Behavior}

\begin{instruction}[Help-Seeking Behavior]
Record each question asked by the four students. 
Write down the question immediately, code the question if possible. 

Use the following codes:
\begin{itemize}
    \item \textbf{Ex}: Question relates to an exercise.
    \item \textbf{Th}: Question relates to theory or concepts.
    \item \textbf{WP}: Question relates to Waterproof itself (code a question about an exercise where the student is using Waterproof as both \textbf{Ex} and \textbf{WP}, questions about technicalities of Waterproof as \textbf{WP}).
    \item \textbf{V}: Question is mostly about validation (e.g. "Is this proof correct?").
    \item \textbf{BTX}: Where X is a number from 1 to 6, corresponding to the number of Bloom's taxonomy that matches the question best.
    \item \textbf{OTM}: Question is off-topic to this course, but related to mathematics.
    \item \textbf{OT}: Question is off-topic, not related to mathematics.
    \item \textbf{U}: TA is unsure of correct code for this question, marks this to be checked later.
\end{itemize}

Record the observations in an Excel sheet in OneDrive with the name questions\_groupnumber\_YYYYMMDD.xlsx.
Enter the codes as a comma-separated list in the Codes column.
The columns should be: Student name, Time (rounded to nearest minute), Question, Codes, Comment.
\end{instruction}

\subsection{General observations}

\begin{instruction}[General Observations]
Record general observations about student behavior, their interactions with exercises with or without Waterproof, as well as any feedback that students give you.
Record the Student Name only if the observation is feedback given by a student.
Record the observations in an Excel sheet in OneDrive with the name general\_teachername\_YYYYMMDD.xlsx.
The columns should be: Time (rounded to nearest minute), Description, Student Name, Comment.

These observations can be more qualitative in nature, and are more of a catch-all.
Feel free to include observations about anything you find interesting. 

\end{instruction}

\section{Entry Survey}\label{app:entry-survey}

\subsection*{About Me}

\begin{enumerate}
    \item\label{q:entry:consent} Did you fill in the consent form?
    \begin{itemize}
        \item[$\circ$] No
        \item[$\circ$] Yes
    \end{itemize}
    \textit{(Students who answered ``No'' were directed to the end of the survey.)}

    \item\label{q:entry:name} What is your name?

    \item\label{q:entry:studentnumber} What is your student number?

    \item\label{q:entry:gender} What is your gender?
    \begin{itemize}
        \item[$\circ$] Male
        \item[$\circ$] Female
        \item[$\circ$] Other, feel free to specify, if you want: \underline{\hspace{4cm}}
        \item[$\circ$] Prefer not to say
    \end{itemize}

    \item\label{q:entry:age} What is your age?

    \item\label{q:entry:programme} What is your study programme?
    \begin{itemize}
        \item[$\circ$] Mathematics (no double bachelor)
        \item[$\circ$] Mathematics and applications
        \item[$\circ$] Mathematics/Physics
        \item[$\circ$] Mathematics/Computer Science
        \item[$\circ$] Mathematics/Economics
        \item[$\circ$] Other: \underline{\hspace{4cm}}
    \end{itemize}

    \item\label{q:entry:years} I have been in higher education for \ldots\ years. (So 0 if you are a first-year student.)
\end{enumerate}

\subsection*{Previous Grades}

We will ask for your grades for central exams and total final grades for secondary education. If you were not in the Dutch secondary system, please enter grades on a 10 point scale that best correspond to the respective grade.

\begin{enumerate}[resume]
    \item\label{q:entry:mathB-central} What was your grade for the central exam for Mathematics~B (or other obligatory mathematics course)?

    \item\label{q:entry:mathB-final} What was your final grade for Mathematics~B (or other obligatory mathematics course)?

    \item\label{q:entry:mathD} What was your final grade for Mathematics~D (or other optional mathematics course)?
    \begin{itemize}
        \item[$\circ$] Grade: \underline{\hspace{3cm}}
        \item[$\circ$] I did not take Mathematics~D, but my school did offer it.
        \item[$\circ$] My school did not offer Mathematics~D.
    \end{itemize}

    \item\label{q:entry:proofs} Have you ever written a mathematical proof, outside the secondary school curriculum? (Select all that apply.)
    \begin{itemize}
        \item[$\square$] No
        \item[$\square$] Yes, in (training for) mathematics competitions
        \item[$\square$] Yes, during a prior study
        \item[$\square$] Yes, during self-study
        \item[$\square$] Yes, in another way: \underline{\hspace{4cm}}
    \end{itemize}
\end{enumerate}

\section{Final Survey}\label{app:final-survey}

\subsection*{About Me}

\begin{enumerate}
    \item\label{q:final:consent} Did you fill in the consent form at the start of the course?
    \begin{itemize}
        \item[$\circ$] No
        \item[$\circ$] Yes
    \end{itemize}
    \textit{(Students who answered ``No'' were directed to the end of the survey.)}

    \item\label{q:final:name} What is your name?

    \item\label{q:final:studentnumber} What is your student number?
\end{enumerate}

\subsection*{Feedback}

\begin{enumerate}[resume]
    \item\label{q:final:feedback} For each of the following statements, choose the most appropriate option. Feedback refers to all kinds of feedback received during the course: from teachers, TAs, written, oral and through tools. Options: Strongly disagree -- Somewhat disagree -- Neither agree nor disagree -- Somewhat agree -- Strongly agree.
    \begin{enumerate}[label=(\alph*)]
        \item\label{q:final:feedback:quantity1} I would have liked to have gotten more feedback during the course.
        \item\label{q:final:feedback:availability} The teacher or teacher assistant were available when I wanted to ask a question.
        \item\label{q:final:feedback:quantity2} I feel like I got a lot of feedback during this course.
        \item\label{q:final:feedback:frequency1} I would have liked to have gotten feedback more frequently during the course.
        \item\label{q:final:feedback:timing} I got feedback exactly when I needed it.
        \item\label{q:final:feedback:frequency2} I generally got feedback at a frequency I found helpful to learning.
        \item\label{q:final:feedback:wanting} I sometimes found myself wanting feedback and not getting it.
        \item\label{q:final:feedback:often} I got feedback quite often.
        \item\label{q:final:feedback:quantity3} I feel like I got a lot of feedback during this course.
    \end{enumerate}
\end{enumerate}

\subsection*{Time}

\begin{enumerate}[resume]
    \item\label{q:final:hours-excl} I spent \ldots\ hours on Bewijzen in de Wiskunde (excluding the 8 hours of weekly lectures and tutorials).

    \item\label{q:final:hours-total} I spent \ldots\ hours on Bewijzen in de Wiskunde in total.

    \item\label{q:final:hours-studies} I spent \ldots\ hours on my studies in total from 1 September 2025 up to and including 4 November 2025.

    \item\label{q:final:exercises} I seriously attempted or completed \ldots\ number of exercises listed in the course planning.
\end{enumerate}

\subsection*{Waterproof in Exercises}

\begin{enumerate}[resume]
    \item\label{q:final:wp-used} I used Waterproof to do 1 or more exercises.
    \begin{itemize}
        \item[$\circ$] No
        \item[$\circ$] Yes
    \end{itemize}
    \textit{(Students who answered ``No'' were directed to the end of this block.)}

    \item\label{q:final:wp-exercises} I seriously attempted or completed \ldots\ exercises using Waterproof.

    \item\label{q:final:wp-likert} Select the option that best matches your opinion. Options: Strongly disagree -- Somewhat disagree -- Neither agree nor disagree -- Somewhat agree -- Strongly agree.
    \begin{enumerate}[label=(\alph*)]
        \item I enjoyed using Waterproof.
        \item Doing exercises in Waterproof helped me develop my proof writing skills.
        \item Using Waterproof was a pleasant experience.
        \item Doing an exercise in Waterproof is easier for me than doing the same exercise on paper (explain why).
    \end{enumerate}

    \item\label{q:final:wp-comments} List three comments about how doing exercises in Waterproof did or did not influence your learning process.

    \item\label{q:final:wp-reasons} List three reasons you continued or stopped using Waterproof.

    \item\label{q:final:wp-other} If you have any other comments regarding the use of Waterproof for exercises, leave them here.
\end{enumerate}

\subsection*{Waterproof in Class}

\begin{enumerate}[resume]
    \item\label{q:final:wp-class-group} I was in a group where Waterproof was mentioned in class (Groups 1 and 2, taught by Johan and Pim).
    \begin{itemize}
        \item[$\circ$] No
        \item[$\circ$] Yes
    \end{itemize}
    \textit{(Students who answered ``No'' were directed to the end of this block.)}

    \item\label{q:final:wp-class-likert} Select the most appropriate option. Options: Strongly disagree -- Somewhat disagree -- Neither agree nor disagree -- Somewhat agree -- Strongly agree.
    \begin{enumerate}[label=(\alph*)]
        \item I would have liked to see more use of Waterproof by the teacher.
        \item The use of Waterproof by the teacher helped me better understand proofs.
        \item The use of Waterproof by the teacher had a negative effect on my learning.
    \end{enumerate}

    \item\label{q:final:wp-class-comments} Give three comments about how the use of Waterproof in class contributed to your learning or not.

    \item\label{q:final:wp-suggestions} If you have any suggestions for the use of Waterproof in this course, or the Waterproof software in general, leave them here.
\end{enumerate}

\subsection*{General}

\begin{enumerate}[resume]
    \item\label{q:final:general} If you have any other remarks concerning the course Bewijzen in de Wiskunde, or anything related to this research, leave them here.
\end{enumerate}

\section{Translation of Dutch proof}\label{app:translation}

\begin{table}[h]
\centering
\begin{tabular}{|l|}
\hline
  We need to prove that $\forall y \in f(A), \exists x \in A, y = g(x)$. \\
  Take $y \in f(A)$. \\
  $f(A) = \{f(x)\mid x \in A\}$, so $\exists x' \in A, y = f(x)$. \\
  Obtain this $x' \in A$. \\
  So it holds that $y \in f(x')$. \\
  We need to prove that $\exists x \in A, y = g(x)$. \\
  Choose $x := x'$. \\
  (Indeed, $x \in A$.) \\
  We need to prove that $y = g(x')$. \\
  By definition of $g$ it holds that $g(x')=f(x')$, so it is true that $y = g(x')$. \\
  We conclude that $g$ is surjective. \\
\hline

\end{tabular}
\caption{Translation of the Dutch proof in \Cref{tab:dutchwp}}
\label{tab:dutchwp_translation}
\end{table}

\end{document}